\input amstex
\documentstyle{amsppt}
\hsize160mm \vsize220mm

 \centerline{\bf On some intermediate mean values}

\vskip 1cm

\centerline{ \bf Slavko Simic}
\bigskip
\centerline{ Mathematical Institute SANU, Kneza Mihaila 36, 11000 Belgrade, Serbia}
\bigskip
\centerline{E-mail: ssimic\@turing.mi.sanu.ac.rs}

 \vskip 1cm


\vskip 1cm

 {\sevenrm{\bf Abstract} We give a necessary and sufficient mean condition
for the quotient of two Jensen functionals and define a new class
$\Lambda_{f,g}(a, b)$ of mean values where $f, g$ are continuously
differentiable convex functions satisfying the relation $f''(t)=t
g''(t), t\in \Bbb R^+$. Then we asked for a characterization of
$f, g$ such that the inequalities $H(a, b)\le \Lambda_{f, g}(a,
b)\le A(a, b)$ or $L(a, b)\le \Lambda_{f, g}(a, b)\le I(a, b)$
hold for each positive $a, b$, where $H, A, L, I$ are the
harmonic, arithmetic, logarithmic and identric means,
respectively. For a subclass of $\Lambda$ with $g''(t)=t^s, \ s\in
\Bbb R$, this problem is thoroughly solved.}
\bigskip
{\bf 1. Introduction}
\bigskip
{\bf 1. 1} \ It is said that the mean $P$ is intermediate relating to the means $M$ and $N, \ M\le N$ if the
relation
$$
M(a, b)\le P(a, b)\le N(a, b),
$$
holds for each two positive numbers $a, b$.
\bigskip
It is also well known that
$$
\min\{a, b\}\le H(a, b)\le G(a, b)\le L(a, b)\le I(a, b)\le A(a,
b)\le S(a,b)\le\max\{a, b\},\eqno (1)
$$
where
$$
H=H(a, b):=2(1/a+1/b)^{-1}; \ \ G=G(a, b):=\sqrt {ab}; \ \ L=L(a,
b):={b-a\over \log b-\log a};
$$
$$
 I=I(a, b):=(b^b/a^a)^{1/(b-a)}/e; \ \ A=A(a, b):={a+b\over 2}; \ \
 S=S(a,b):= a^{a\over a+b} b^{b\over a+b},
$$
are the harmonic, geometric, logarithmic, identric, arithmetic and
Gini  mean, respectively.
\bigskip
An easy task is to construct intermediate means related to two
given means $M$ and $N$ with $M\le N$. For instance, for an
arbitrary mean $P$, we have that
$$
M(a, b)\le P(M(a,b),N(a,b))\le N(a, b).
$$
\bigskip
The problem is more difficult if we have to decide whether the given mean is intermediate or not. For example,
the relation
$$
L(a, b)\le S_s(a, b)\le I(a, b),
$$
holds for each positive $a$ and $b$ if and only if $0\le s\le 1$,
where the Stolarsky mean $S_s$ is defined by (cf [4])
$$
S_s(a,b):=\Bigl({b^s-a^s\over s(b-a)}\Bigr)^{1/(s-1)}.
$$
\bigskip
Also,
$$
G(a, b)\le A_s(a, b)\le A(a, b),
$$
holds if and only if $0\le s\le 1$, where the H\"{o}lder mean of order $s$ is defined by
$$
A_s(a, b):=\Bigl({a^s+b^s\over 2}\Bigr)^{1/s}.
$$
\bigskip
An inverse problem is to find best possible approximation of a
given mean $P$ by elements of an ordered class of means $S$. A
good example for this topic is comparison between the logarithmic
mean and the class $A_s$ of H\"{o}lder means of order $s$. Namely,
since $A_0=\lim_{s\to 0}A_s=G$ and $A_1=A$, it follows from $(1)$
that
$$
A_0\le L\le A_1.
$$
\bigskip
Since $A_s$ is monotone increasing in $s$, an improving of the
above is given by Carlson [2]:
$$
A_0\le L\le A_{1/2}.
$$
\bigskip
Finally,  Lin shoved in [3] that
$$
A_0\le L\le A_{1/3},
$$
is the best possible approximation of the logarithmic mean by the
means from the class $A_s$.
\bigskip

Numerous similar results have been obtained recently. For example,
an approximation of Seiffert's mean  by the class $A_s$ is given
in [6], [8].

\bigskip

In this article we shall give best possible approximations for a
whole variety of elementary means $(1)$ by the class $\lambda_s$
defined below (see Thm 3.).

\bigskip

{\bf 1. 2.} \ Let $f, g$ be twice continuously differentiable
(strictly) convex functions on $\Bbb R^+$. By definition (cf [1],
p. 5),
$$
\bar f(a, b):=f(a)+f(b)-2f({a+b\over 2})>0, \ a\neq b,
$$
and
$$
\bar f(a, b)=0,
$$
if and only if $a=b$.
\bigskip
It turns out that the expression
$$
\Lambda_{f, g}(a, b):={\bar f(a, b)\over \bar g(a, b)}={f(a)+f(b)-2f({a+b\over 2})\over g(a)+g(b)-2g({a+b\over
2})},
$$
represents a mean of two positive numbers $a, b$; that is, the relation
$$
\min\{a, b\}\le \Lambda_{f, g}(a, b)\le \max \{a, b\},\eqno (2)
$$
holds for each $a, b\in \Bbb R^+$, if and only if the relation
$$
f''(t)=t g''(t),\eqno (3)
$$
holds for each $t\in \Bbb R^+$.
\bigskip
Let $f, g\in C^\infty(0,\infty)$ and denote by $\Lambda$ the set $\{(f, g)\}$ of convex functions satisfying the
relation (3). There is a natural question how to improve the bounds in (2); in this sense we come upon the
following intermediate mean problem:
\bigskip
{\bf Open question} \ {\it Under what additional conditions on $f, g\in\Lambda$, the inequalities
$$
H(a, b)\le \Lambda_{f, g}(a, b)\le A(a, b),
$$
or, more tightly,
$$
L(a, b)\le \Lambda_{f, g}(a, b)\le I(a, b),
$$
hold for each $a, b\in \Bbb R^+$?}
\bigskip
As an illustration, consider the function $f_s(t)$ defined to be
$$
f_s(t)=\cases {(t^s-st+s-1)/s(s-1)} &, s(s-1)\neq 0;\\
              t-\log t-1 &, s=0;\\
              t\log t-t+1 &, s=1.
              \endcases
              $$
              \bigskip
Since
$$
f'_s(t)=\cases {t^{s-1}-1\over s-1} &, s(s-1)\neq 0;\\
               1-{1\over t} &, s=0;\\
               \log t &, s=1,
               \endcases
               $$
and
$$
f''_s(t)=t^{s-2}, \ s\in \Bbb R, \ t>0,
$$
it follows that $f_s(t)$ is a twice continuously differentiable convex function for $s\in \Bbb R, \ t\in \Bbb
R^+$.
\bigskip
Moreover, it is evident that $(f_{s+1}, f_s)\in \Lambda$.
\bigskip
 We shall give in the sequel a complete answer to the above question concerning the means
 $$
 \bar f_{s+1}(a, b)/\bar f_s(a, b):=\lambda_s(a, b)
 $$
  defined by
$$
\lambda_s(a, b)=\cases {s-1\over s+1}{a^{s+1}+b^{s+1}-2({a+b\over 2})^{s+1}\over a^s+b^s-2({a+b\over 2})^s},&
s\in \Bbb
R/\{-1, 0, 1\};\\
                {2\log {a+b\over 2}-\log a-\log b\over {1\over 2a}+{1\over 2b}-{2\over a+b}}, & s=-1;\\
                {a\log a+b\log b-(a+b)\log{a+b\over 2}\over 2\log{a+b\over 2}-\log a-\log b},& s=0;\\
                {(b-a)^2\over 4(a\log a+b\log b-(a+b)\log{a+b\over 2})},& s=1.
                \endcases
                $$
\bigskip

Those means are obviously symmetric and homogeneous of order one.

\bigskip

 As a consequence we obtain some new intermediate mean
values; for instance, we show that the inequalities
$$
H(a, b)\le \lambda_{-1}(a, b)\le G(a, b)\le \lambda_0(a, b)\le L(a,b)\le \lambda_1(a, b)\le I(a, b),
$$
hold for arbitrary $a, b\in \Bbb R^+$.
\bigskip
Note that

$$
\lambda_{-1}={2G^2\log(A/G)\over A-H}; \
\lambda_0=A{\log(S/A)\over \log(A/G)}; \ \lambda_1={1\over
2}{A-H\over \log(S/A)}.
$$

\bigskip

{\bf 2. Results}
\bigskip
We prove firstly the following
\bigskip
{\bf Theorem 1} \ {\it Let $f, g\in C^2(I)$ with $g''>0$. The
expression $\Lambda_{f, g}(a, b)$ represents a mean of arbitrary
numbers $a, b\in I$ if and only if the relation
$$
f''(t)=t g''(t)\eqno (3)
$$
holds for $t\in I$. }
\bigskip
{\bf Remark 1} \ {\it In the same way, for arbitrary $p,q>0, p+q=1$, it can be deduced that the quotient
$$
\Lambda_{f, g}(p,q; a, b):={pf(a)+qf(b)-f(pa+qb)\over pg(a)+qg(b)-g(pa+qb)}
$$
represents a mean value of numbers $a, b$ if and only if (3)
holds.}
\bigskip
A generalization of the above assertion is the next
\bigskip
 {\bf Theorem 2} \ {\it Let $f, g: I\to \Bbb R$ be twice continuously differentiable functions with
 $g''>0$ on $I$ and let $p=\{p_i\}, i=1, 2, \cdots, \ \sum p_i=1$ be an arbitrary
positive weight sequence. Then the quotient of two Jensen
functionals
$$
\Lambda_{f, g}(p, x):={\sum_1^n p_i f(x_i)-f(\sum_1^n p_i x_i)\over \sum_1^n p_i g(x_i)-g(\sum_1^n p_i x_i)}, \
\ n\ge 2,
$$
represents a mean of an arbitrary set of real numbers $x_1, x_2,
\cdots, x_n\in I$ if and only if the relation
$$
f''(t)=tg''(t)
$$
holds for each $t\in I$.}
\bigskip
{\bf Remark 2} \ {\it It should be noted that the relation
$f''(t)=tg''(t)$ determines $f$ in terms of $g$ in an easy way.
Precisely,

$$
f(t)=tg(t)-2G(t)+ct+d,
$$

where $G(t):=\int_1^t g(u)du$ and $c$ and $d$ are constants.}
\bigskip
Our results concerning the means $\lambda_s(a,b), \ s\in \Bbb R$ are included in the following
\bigskip
{\bf Theorem 3} \ {\it For the class of means $\lambda_s(a, b)$
defined above, the following assertions hold for each $a, b\in
\Bbb R^+$.
\bigskip
{\bf 1.} \ \ The means $\lambda_s(a, b)$ are monotone increasing in $s$;
\bigskip
{\bf 2.} \ $\lambda_s(a, b)\le H(a, b)$ for each $s\le -4$;
\bigskip
{\bf 3.} \ $H(a, b)\le \lambda_s(a, b)\le G(a, b)$ for $-3\le s\le -1$;
\bigskip
{\bf 4.} \ $G(a, b)\le \lambda_s(a, b)\le L(a, b)$ for $-1/2\le s\le 0$;
\bigskip
{\bf 5.} \ there is a number $s_0\in (1/12, 1/11)$ such that \ $L(a, b)\le \lambda_s(a, b)\le I(a,b)$ for
$s_0\le s\le 1$;
\bigskip
{\bf 6.} \ there is a number $s_1\in (1.03, 1.04)$ such that \ $I(a, b)\le \lambda_s(a, b)\le A(a, b)$ for
$s_1\le s\le 2$;
\bigskip
{\bf 7.} \ $A(a, b)\le \lambda_s(a, b)\le S(a,b)$ for each $2\le
s\le 5$;
\bigskip
{\bf 8.} \ there is no finite $s$ such that the inequality
$S(a,b)\le \lambda_s(a,b)$ holds for each $a,b\in \Bbb R^+$.
\bigskip
The above estimations are best possible. }
\bigskip
{\bf 3. Proofs}
\bigskip
{\bf Proof of Theorem 1} \ We prove firstly the necessity of the condition (3).
\bigskip
Since $\Lambda_{f, g}(a, b)$ is a mean value for arbitrary $a,
b\in I; \ a\neq b$, we have
$$
\min\{a, b\}\le \Lambda_{f, g}(a, b)\le \max\{a, b\}.
$$
Hence
$$
\lim_{b\to a}\Lambda_{f, g}(a, b)=a. \eqno (4)
$$
\bigskip
From the other hand, due to l'Hospital's rule we obtain
$$
\lim_{b\to a}\Lambda_{f, g}(a, b)=\lim_{b\to a}\Bigl({f'(b)-f'({a+b\over 2})\over g'(b)-g'({a+b\over 2})}\Bigr)
=\lim_{b\to a}\Bigl({2f''(b)-f''({a+b\over 2})\over 2g''(b)-g''({a+b\over 2})}\Bigr)
$$
$$
={f''(a)\over g''(a)}. \eqno (5)
$$
Comparing (4) and (5) the desired result follows.
\bigskip
Suppose now that (3) holds and let $a<b$. Since $g''(t)>0 \ t\in
[a, b]$ by {\it the Cauchy mean value theorem} there exists
$\xi\in ({a+t\over 2}, t)$ such that
$$
{f'(t)-f'({a+t\over 2})\over g'(t)-g'({a+t\over 2})}={f''(\xi)\over g''(\xi)}=\xi. \eqno (6)
$$
\bigskip
But,
$$
a\le {a+t\over 2}<\xi<t\le b,
$$
and, since $g'$ is strictly increasing, \ $g'(t)-g'({a+t\over 2})>
0, \ t\in [a, b]$.
\bigskip
Therefore, by (6) we get
$$
a( g'(t)-g'({a+t\over 2}))\le  f'(t)-f'({a+t\over 2})\le b( g'(t)-g'({a+t\over 2})). \eqno (7)
$$
\bigskip
Finally, integrating (7) over $t\in [a, b]$ we obtain the assertion from Theorem 1.
\bigskip

 {\bf Proof of Theorem 2} \ We shall give a proof of this assertion by induction on $n$.
\bigskip
By Remark 1, it holds for $n=2$.
\bigskip
Next, it is not difficult to check the identity
$$
\sum_1^n p_i f(x_i)-f(\sum_1^n p_i x_i)=(1-p_n)(\sum_1^{n-1} p'_i f(x_i)-f(\sum_1^{n-1} p'_i
x_i))
$$
$$
+[(1-p_n)f(T)+p_nf(x_n)-f((1-p_n)T+p_n x_n)],
$$
where
$$
T:=\sum_1^{n-1}p_i'x_i; \ \ p_i':=p_i/(1-p_n), \ \ i=1, 2, \cdots, n-1; \ \sum_1^{n-1}p'_i=1.
$$
\bigskip
Therefore, by induction hypothesis and Remark 1, we get
$$
\sum_1^n p_i f(x_i)-f(\sum_1^n p_i x_i)\le \max\{x_1, x_2,\cdots x_{n-1}\}(1-p_n)(\sum_1^{n-1} p'_i
g(x_i)-g(\sum_1^{n-1} p'_i x_i))
$$
$$
+\max\{T, x_n\}[(1-p_n)g(T)+p_ng(x_n)-g((1-p_n)T+p_n x_n)]
$$
$$
\le \max \{x_1, x_2,\cdots, x_n\}((1-p_n)(\sum_1^{n-1} p'_i g(x_i)-g(\sum_1^{n-1} p'_i x_i))
$$
$$
+[(1-p_n)g(T)+p_ng(x_n)-g((1-p_n)T+p_n x_n)])
$$
$$
=\max \{x_1, x_2,\cdots, x_n\}(\sum_1^n p_i g(x_i)-g(\sum_1^n p_i x_i)).
$$
\bigskip
The inequality
$$
\min\{x_1, x_2, \cdots, x_n\}\le \Lambda_{f, g}(p, x),
$$
can be proved analogously.

\bigskip

For the proof of necessity, put $x_2=x_3=\dots =x_n$ and proceed
as in Theorem 1.

\bigskip

{\bf Remark} \ {\it It is evident from $(3)$ that if $I\subseteq
\Bbb R^+$ then $f$ has to be also convex on $I$. Otherwise, it
shouldn't be the case. For example, the conditions of Theorem 2
are satisfied with $f(t)=t^3/3, g(t)=t^2, t\in \Bbb R$. Hence, for
an arbitrary sequence $\{x_i\}_1^n$ of real numbers, we obtain
$$
\min\{x_1, x_2, \cdots, x_n\}\le {\sum_1^n p_i x_i^3-(\sum_1^n p_i
x_i)^3\over 3(\sum_1^n p_i x_i^2-(\sum_1^n p_i x_i)^2)} \le
\max\{x_1, x_2, \cdots, x_n\}.
$$}

 Because the above inequality does not depend on $n$, a probabilistic interpretation of the above result is
 contained in the following
\bigskip
{\bf Theorem 4.} \ {\it For an arbitrary probability law $F$ of
random variable $X$ with support on $(-\infty, +\infty)$, we have
$$
(EX)^3+3(\min X) \ \sigma_X^2\le EX^3\le (EX)^3+3 (\max X) \
\sigma_X^2.
$$}

\bigskip

 {\bf Proof of Theorem 3, part 1} \ We shall prove a general assertion of this type. Namely, for an
arbitrary positive sequence $\bold x=\{x_i\}$ and an associated weight sequence $\bold p=\{p_i\}, \ i=1, 2,
\cdots$, denote
$$
\chi_s(\bold p, \bold x):=\cases {\sum p_i x_i^s-(\sum p_i x_i)^s\over s(s-1)}, & s\in \Bbb R/\{0,
1\};\\
                                     \log(\sum p_i x_i)-\sum p_i\log x_i, & s=0;\\
                                     \sum p_i x_i\log x_i-(\sum p_i x_i)\log (\sum p_i x_i), & s=1.
                                     \endcases
                                     $$
\bigskip
For $s\in \Bbb R, r>0$ we have
$$
\chi_s(\bold p, \bold x) \chi_{s+r+1}(\bold p, \bold x)\ge \chi_{s+1}(\bold p, \bold x)\chi_{s+r}(\bold p, \bold
x),\eqno (4)
$$
which is equivalent to
\bigskip
{\bf Theorem 3a} \ {\it The sequence $\{\chi_{s+1}(\bold p, \bold x)/\chi_s(\bold p, \bold x)\}$ is monotone
increasing in $s, \ s\in \Bbb R$.}
\bigskip
 This assertion follows applying the result from ([5], Theorem 2) which states that
 \bigskip
{\bf Lemma 1} \ {\it For $-\infty<a<b<c<+\infty$, the inequality
$$
(\chi_b(\bold p, \bold x))^{c-a}\le (\chi_a(\bold p, \bold x))^{c-b}(\chi_c(\bold p, \bold x))^{b-a},
$$
holds for arbitrary sequences $\bold p, \bold x$.}
\bigskip
Putting there $a=s, b=s+1, c=s+r+1$ and $a=s, b=s+r, c=s+r+1$, we successively  obtain
$$
(\chi_{s+1}(\bold p, \bold x))^{r+1}\le (\chi_s(\bold p, \bold x))^r \chi_{s+r+1}(\bold p, \bold x),
$$
and
$$
(\chi_{s+r}(\bold p, \bold x))^{r+1}\le \chi_s(\bold p, \bold x)(\chi_{s+r+1}(\bold p, \bold x))^r.
$$
\bigskip
Since $r>0$, multiplying those inequalities we get the relation (4) i. e. the proof of Theorem 3a.
\bigskip
The part 1. of Theorem 3 follows for $p_1=p_2=1/2$.
\bigskip
A general way to prove the rest of Theorem 3 is to use an easy-checkable identity
$$
{\lambda_s(a, b)\over A(a, b)}=\lambda_s(1+t, 1-t),
$$
with $t:={b-a\over b+a}$.
\bigskip
Since $0<a<b$, we get $0<t<1$. Also,
$$
{H(a, b)\over A(a, b)}=1-t^2; \ \ {G(a, b)\over A(a, b)}=\sqrt
{1-t^2}; \ \ {L(a, b)\over A(a, b)}={2t\over \log(1+t)-\log(1-t)};
\eqno(5)
$$
$$
{I(a, b)\over A(a,b)}=\exp({(1+t)\log(1+t)-(1-t)\log(1-t)\over
2t}-1); \ {S(a, b)\over A(a,b)}=\exp({1\over
2}((1+t)\log(1+t)+(1-t)\log(1-t))).
$$
\bigskip
Therefore, we have to compare some one-variable inequalities and to check their validness for each $t\in (0,
1)$.
\bigskip
For example, we shall prove that the inequality
$$
\lambda_s(a, b)\le L(a,b)
$$
holds for each positive $a, b$ if and only if $s\le 0$.
\bigskip
Since $\lambda_s(a, b)$ is monotone increasing in $s$, it is enough to prove that
$$
{\lambda_0(a,b)\over L(a, b)}\le 1.
$$
\bigskip
By the above formulae, this is equivalent to the assertion that the inequality
$$
\phi(t)\le 0\eqno (6)
$$
holds for each $t\in (0, 1)$, with

$$
\phi (t):={\log(1+t)-\log(1-t)\over 2t}((1+t)\log(1+t)+(1-t)\log(1-t))+\log(1+t)+\log(1-t).
$$
\bigskip
We shall prove that the power series expansion of $\phi(t)$ have
non-positive coefficients. Thus the relation (6) will be proved.
\bigskip
Since
$$
{\log(1+t)-\log(1-t)\over 2t}=\sum_0^\infty {t^{2k}\over 2k+1}; \ \ \ \log(1+t)+\log(1-t)=-t^2 \sum_0^\infty
{t^{2k}\over k+1};
$$
$$
(1+t)\log(1+t)+(1-t)\log(1-t)=t^2 \sum_0^\infty {t^{2k}\over (k+1)(2k+1)},
$$
we get
$$
\phi(t)/t^2=\sum_{n=0}^\infty(-{1\over n+1}+\sum_{k=0}^n {1\over (2n-2k+1)(k+1)(2k+1)})t^{2n}=\sum_0^\infty c_n
t^{2n}.
$$
\bigskip
Hence,
$$
c_0=c_1=0; \ \ c_2=-1/90,
$$
and, after some calculation, we get
$$
c_n={2\over (n+1)(2n+3)}\Bigl((n+2)\sum_1^n{1\over 2k+1}-(n+1)\sum_1^n{1\over 2k}\Bigr), \ n>1.
$$
\bigskip
Now, one can easily prove (by induction, for example) that
$$
d_n:=(n+2)\sum_1^n{1\over 2k+1}-(n+1)\sum_1^n{1\over 2k},
$$
is a negative real number for $n\ge 2$. Therefore $c_n\le 0$, and the proof of the first part is done.
\bigskip
For $0<s<1$ we have
$$
{\lambda_s(a, b)\over L(a, b)}-1=\frac{(1-s)((1+t)^{s+1}+(1-t)^{s+1}-2)\log
\frac{1+t}{1-t}}{2t(1+s)(2-(1+t)^s-(1-t)^s)}-1={1\over 6}st^2+O(t^4) \ \ \ (t\to 0).
$$
\bigskip
Therefore, $\lambda_s(a, b)>L(a, b)$ \ for $s>0$ and sufficiently small $t:=(b-a)/(b+a)$.
\bigskip
Similarly, we shall prove that the inequality
$$
\lambda_s(a, b)\le I(a, b),
$$
holds for each $a, b; 0<a<b$ if and only if $s\le 1$.
\bigskip
As before, it is enough to consider the expression
$$
{I(a, b)\over \lambda_1(a, b)}=e^{\mu(t)}\nu(t):=\psi(t),
$$
with
$$
\mu(t)={(1+t)\log(1+t)-(1-t)\log(1-t)\over 2t}-1; \ \ \nu(t)={(1+t)\log(1+t)+(1-t)\log(1-t)\over t^2}.
$$
\bigskip
It is not difficult to check the identity
$$
\psi'(t)=-e^{\mu(t)}\phi(t)/t^3.
$$
\bigskip
Hence by (6), we get $\psi'(t)>0$ i. e. $\psi(t)$ is monotone increasing for $t\in(0, 1)$.
\bigskip
Therefore
$$
{I(a, b)\over \lambda_1(a, b)}\ge \lim_{t\to 0^+}\psi(t)=1.
$$
\bigskip
By monotonicity it follows that $\lambda_s(a, b)\le I(a, b)$ for $s\le 1$.
\bigskip
For $s>1, \ \ {b-a\over b+a}=t$, we have
$$
\lambda_s(a,b)-I(a, b)=({1\over 6}(s-1)t^2+O(t^4))A(a,b) \ \ \
(t\to 0^+).
$$
\bigskip
Hence, $\lambda_s(a,b)>I(a, b)$ for $s>1$ and $t$ sufficiently small .
\bigskip
From the other hand,
$$
\lim_{t\to 1^-}[{\lambda_s(a, b)\over I(a, b)}-1]={e(s-1)(2^{s+1}-2)\over 2(s+1)(2^s-2)}-1:=\tau(s).
$$
\bigskip
Examining the function $\tau(s)$, we find out that it has the only real zero at $s_0\approx 1.0376$ and is
negative for $s\in(1, s_0)$.
\bigskip
 {\bf Remark 2} \ {\it Since $\psi(t)$ is monotone increasing, we also get
$$
{I(a, b)\over \lambda_1(a, b)}\le \lim_{t\to 1^-}\psi(t)={4\log 2\over e}.
$$
Hence
$$
1\le {I(a, b)\over \lambda_1(a, b)}\le {4\log 2\over e}.
$$

A calculation gives ${4\log 2\over e}\approx 1.0200$.}

\bigskip

Note also that
$$
\lambda_2(a, b)\equiv A(a, b).
$$
\bigskip
Therefore, applying the assertion from the part 1., we get
$$
\lambda_s(a, b)\le A(a, b), \ s\le 2; \ \ \lambda_s(a, b)\ge A(a, b), \ s\ge 2.
$$

\bigskip

Finally, we give a detailed proof of the part 7.

\vskip 1cm

We have to prove that $\lambda _s(a,b)\le S(a,b)$ for $s\le 5$.
Since $\lambda_s(a,b)$ is monotone increasing in $s$, it is
sufficient to prove that the inequality
$$
\lambda _5(a,b)\le S(a,b)
$$

holds for each $a,b\in \Bbb R^+$.

\bigskip

Therefore, by the transformation given above, we get

$$
\log {\lambda_5\over A}=\log\bigl[{2\over
3}{(1+t)^6+(1-t)^6-2\over
(1+t)^5+(1-t)^5-2}\bigr]=\log\bigl[{2\over 15}{15+15t^2+t^4\over
2+t^2}\bigr]
$$
$$
\le \log\bigl[{1+t^2+t^4/4\over
1+t^2/2}\bigr]=\log(1+t^2/2)=t^2/2-t^4/8+t^6/24-\cdots
$$
$$
\le t^2/2+t^4/12+t^6/30+\cdots={1\over
2}((1+t)\log(1+t)+(1-t)\log(1-t))=\log{S/A},
$$

and the proof is done.

\bigskip

Further, we have to show that $\lambda_s(a,b)>S(a,b)$ for some
positive $a,b$ whenever $s>5$.

\bigskip

Indeed, since
$$
(1+t)^s+(1-t)^s-2={s\choose 2}t^2+{s\choose 4}t^4+O(t^6),
$$

for $s>5$ and sufficiently small $t$, we get

$$
{\lambda_s\over A}={s-1\over s+1}{{s+1\choose 2}t^2+{s+1\choose
4}t^4+O(t^6)\over {s\choose 2}t^2+{s\choose 4}t^4+O(t^6)}
$$
$$
={1+(s-1)(s-2)t^2/12+O(t^4)\over
1+(s-2)(s-3)t^2/12+O(t^4)}=1+({s\over 6}-{1\over 3})t^2+O(t^4).
$$

\bigskip

Similarly,

$$
{S\over A}=\exp({1\over
2}((1+t)\log(1+t)+(1-t)\log(1-t)))=\exp(t^2/2+O(t^4))=1+t^2/2+O(t^4).
$$

Hence,

$$
{1\over A}(\lambda_s-S)={1\over 6}(s-5)t^2+O(t^4),
$$

and this expression is positive for $s>5$ and $t$ sufficiently
small, i.e. $a$ sufficiently close to $b$.

\bigskip

As for the part 8., applying the above transformation we obtain

$$
{\lambda _s(a,b)\over S(a,b)}={s-1\over
s+1}{(1+t)^{s+1}+(1-t)^{s+1}-2\over
(1+t)^s+(1-t)^s-2}\exp(-{1\over
2}((1+t)\log(1+t)+(1-t)\log(1-t))),
$$

where $0<a<b, \ t={b-a\over b+a}$.

\bigskip

Since for $s>5$,

$$
\lim_{t\to 1^-}{\lambda_s\over S}={s-1\over s+1}{2^s-1\over
2^s-2},
$$
and the last expression is less than one, it follows that the
inequality $S(a,b)<\lambda_s(a,b)$ cannot hold whenever ${b\over
a}$ is sufficiently large.

\bigskip

The rest of the proof is straightforward.

\bigskip

{\bf Acknowledgment} \ The author is indebted to the referees for
valuable suggestions.

 \vskip 1cm

{\bf References}
\bigskip
[1] \ Hardy, G.H., Littlewood, J. E., P\"{o}lya, G.: {\it Inequalities}, Camb. Univ. Press, Cambridge (1978).
\bigskip
[2] \ Carlson, B. C. : {\it The logarithmic mean}, Amer. Math.
Monthly, 79 (1972), pp. 615-618.
\bigskip
[3] \ Lin, T. P. : {\it The power mean and the logarithmic mean},
Amer. Math. Monthly, 81 (1974), pp. 879-883.
\bigskip
[4] \  Stolarsky, K.: \ {\it Generalizations of the logarithmic mean}, Math. Mag. 48 (1975), pp. 87-92.
\bigskip
[5] \ Simic, S. : {\it On logarithmic convexity for differences of
power means}, J. Inequal. Appl. Article ID 37359 (2007).
\bigskip
[6] \ Hasto, P.A. : {\it Optimal inequalities between Seiffert's
mean and power means}, Math. Inequal. Appl. Vol. 7, No. 1 (2004),
pp. 47-53.
\bigskip
[7] \ Simic, S. : {\it An extension of Stolarsky means to the
multivariable case}, Int. J. Math. Math. Sci. Article ID 432857
(2009).
\bigskip
[8] \ Yang, Z-H. : {\it Sharp bounds for the second Seiffert mean
in terms of power mean}, arXiv: 1206.5494v1 [math. CA] (2012).

\end